\documentclass{amsart}
\usepackage{amssymb}
\usepackage{amsmath}
\usepackage{amscd}

\newtheorem{prop}{Proposition}[section]

\newtheorem{thm}[prop]{Theorem}

\theoremstyle{definition}

\newcommand{\Aut}{{\mathrm {Aut}}}

\newcommand{\Hom}{{\mathrm {Hom}}}

\newcommand{\Nm}{{\mathrm {Nm}}}

\newcommand{\tr}{{\mathrm {tr}}}

\newcommand{\Gal}{\mathrm {Gal}}

\newcommand{\A}{{\mathbb A}}
\newcommand{\CC}{{\mathbb C}}

\newcommand{\QQ}{{\mathbb Q}}

\newcommand{\ZZ}{{\mathbb Z}}

\newcommand{\TT}{{\mathbb T}}

\newcommand{\WWW}{{\mathcal W}}

\newcommand{\pp}{{\mathfrak p}}

\newcommand{\mm}{{\mathfrak m}}

\newcommand{\FF}{{\mathbb F}}

\newcommand{\GL}{\mathrm {GL}}

\newcommand{\Qbar}{\overline{\mathbb Q}}

\newcommand{\rhobar}{\overline{\rho}}

\DeclareFontEncoding{OT2}{}{} % to enable usage of cyrillic fonts
  \newcommand{\textcyr}[1]{%
    {\fontencoding{OT2}\fontfamily{wncyr}\fontseries{m}\fontshape{n}%
     \selectfont #1}}
\newcommand{\Sha}{{\mbox{\textcyr{Sh}}}}

\begin{document}
\title{Lifting Congruences to Weight $3/2$}
\author{Neil Dummigan}
\author{Srilakshmi Krishnamoorthy}
\date{April 8th, 2015.}

\address{University of Sheffield\\ School of Mathematics and Statistics\\
Hicks Building\\ Hounsfield Road\\ Sheffield, S3 7RH\\
U.K.}
\email{n.p.dummigan@shef.ac.uk}
\address{IIT Madras, Chennai, India.}
\email{srilakshmi@iitm.ac.in}

\begin{abstract}
Given a congruence of Hecke eigenvalues between newforms of weight $2$, we prove, under certain conditions, a congruence between corresponding weight-$3/2$ forms.
\end{abstract}

\maketitle
\section{Introduction}
Let $f=\sum_{n=1}^{\infty}a_n(f)q^n$ and $g=\sum_{n=1}^{\infty}a_n(g)q^n$ be normalised newforms of weight $2$ for $\Gamma_0(N)$, where $N$ is square-free. For each prime $p\mid N$, let $w_p(f)$ and $w_p(g)$ be the eigenvalues of the Atkin-Lehner involution $W_p$ acting on $f$ and $g$, respectively. Write $N=DM$, where $w_p(f)=w_p(g)=-1$ for primes $p\mid D$ and $w_p(f)=w_p(g)=1$ for primes $p\mid M$. We suppose that the number of primes dividing $D$ is odd. (In particular, the signs in the functional equations of $L(f,s)$ and $L(g,s)$ are both $+1$.) Let $B$ be the quaternion algebra over $\QQ$ ramified at $\infty$ and at the primes dividing $D$, with canonical anti-involution $x\mapsto \overline{x}$, $\tr(x):=x+\overline{x}$ and $\Nm(x):=x\overline{x}$. Let $R$ be a fixed Eichler order of level $N$ in a maximal order of $B$. Let $\phi_f, \phi_g$ (determined up to non-zero scalars) be ($\CC$-valued) functions on the finite set $B^{\times}(\QQ)\backslash B^{\times}(\A_f)/\hat{R}$ corresponding to $f$ and $g$ via the Jacquet-Langlands correspondence, where $\A_f$ is the ``finite'' part of the adele ring of $\QQ$ and $\hat{R}=R\otimes_{\ZZ}\hat{\ZZ}$. Let $\{y_i\}_{i=1}^h$ be a set of representatives in $B^{\times}(\A_f)$ of $B^{\times}(\QQ)\backslash B^{\times}(\A_f)/\hat{R}$, $R_i:=B^{\times}(\QQ)\cap(y_i\hat{R}y_i^{-1})$ and $w_i:=|R_i^{\times}|$. Let $L_i:=\{x\in \ZZ+2R_i : \tr(x)=0\}$, and $\theta_i:=\sum_{x\in L_i}q^{\Nm(x)}$, where $q=e^{2\pi iz}$, for $z$ in the complex upper half-plane. For $\phi=\phi_f$ or $\phi_g$, let $$\WWW(\phi):=\sum_{i=1}^h \phi(y_i)\theta_i.$$ This is Waldspurger's theta-lift \cite{W1}, and the Shimura correspondence \cite{Sh} takes $\WWW(\phi_f)$ and $\WWW(\phi_g)$, which are cusp forms of weight $3/2$ for $\Gamma_0(4N)$, to $f$ and $g$, respectively (if $\WWW(\phi_f)$ and $\WWW(\phi_g)$ are non-zero). In the case that $N$ is odd (and square-free), $\WWW(\phi_f)$ and $\WWW(\phi_g)$ are, if non-zero, the unique (up to scaling) elements of Kohnen's space $S_{3/2}^+(\Gamma_0(4N))$ mapping to $f$ and $g$ under the Shimura correspondence \cite{K}. Still in the case that $N$ is odd, $\WWW(\phi_f)\neq 0$ if and only if $L(f,1)\neq 0$, by a theorem of B\"ocherer and Schulze-Pillot \cite[Corollary, p.379]{BS1}, proved by Gross in the case that $N$ is prime \cite[\S 13]{G}.

B\"ocherer and Schulze-Pillot's version of Waldspurger's Theorem \cite{W2},\cite[Theorem 3.2]{BS2} is that for any fundamental discriminant $-d<0$,
$$\sqrt{d}\left(\prod_{p\mid \frac{N}{\mathrm{gcd}(N,d)}}\left(1+\left(\frac{-d}{p}\right )w_p(f)\right )\right)L(f,1)L(f,\chi_{-d},1)=\frac{4\pi^2\langle f,f\rangle}{\langle \phi_f,\phi_f\rangle} (a(\WWW(\phi_f),d))^2,$$
and similarly for $g$, where $\WWW(\phi_f)=\sum_{n=1}^{\infty}a(\WWW(\phi_f),n)q^n$, $\langle f,f\rangle$ is the Petersson norm and $\langle\phi_f,\phi_f\rangle=\sum_{i=1}^h w_i|\phi_f(y_i)|^2$. (They scale $\phi_f$ in such a way that $\langle \phi_f,\phi_f\rangle=1$, so it does not appear in their formula.)

The main goal of this paper is to prove the following.
\begin{thm}\label{main} Let $f,g,\WWW(\phi_f),\WWW(\phi_g)$, $N=DM$ be as above (with $N$ not necessarily odd). Suppose now that $D=q$ is prime. Let $\ell>2$ be a prime such that $\ell\nmid N(q-1)$ or $\ell=q$. Suppose that, for some unramified divisor $\lambda\mid\ell$ in a sufficiently large number field,
$$a_p(f)\equiv a_p(g)\pmod{\lambda}~~~\forall\text{ primes } p\nmid N\ell,$$
and that the residual Galois representation $\rhobar_{f,\lambda}:\Gal(\Qbar/\QQ)\rightarrow\GL_2(\FF_{\lambda})$ is irreducible. Then (with a suitable choice of scaling, such that $\phi_f$ and $\phi_g$ are integral but not divisible by $\lambda$)
$$a(\WWW(\phi_f),n)\equiv a(\WWW(\phi_g),n)\pmod{\lambda}~~~\forall n.$$
\end{thm}
{\em Remarks.}
\begin{enumerate}
\item Though $\phi_f$ and $\phi_g$ are not divisible by $\lambda$, we can still imagine that $\WWW(\phi_f)=\sum_{i=1}^h \phi_f(y_i)\theta_i$ and $\WWW(\phi_g)=\sum_{i=1}^h \phi_g(y_i)\theta_i$ could have all their Fourier coefficients divisible by $\lambda$, so the congruence could be just $0\equiv 0\pmod{\lambda}$ for all $n$. However, unless $\WWW(\phi_f)=\WWW(\phi_g)=0$, this kind of mod $\ell$ linear dependence of the $\theta_i$ seems unlikely, and one might guess that it never happens. This seems related to a conjecture of Kolyvagin, about non-divisibility of orders of Shafarevich-Tate groups of quadratic twists, discussed by Prasanna \cite{P}.
\item The discussion in \cite[\S\S 5.2,5.3]{P} is also relevant to the subject of this paper. In particular, our congruence may be viewed as a square root of a congruence between algebraic parts of $L$-values. Such congruences may be proved in greater generality, as in \cite[Theorem 0.2]{V}, but do not imply ours, since square roots are determined only up to sign. The idea for Theorem \ref{main} came in fact from work of Quattrini \cite[\S 3]{Q}, who proved something similar for congruences between cusp forms and Eisenstein series at prime level, using results of Mazur \cite{M} and Emerton \cite{Em} on the Eisenstein ideal. See Theorem 3.6, and the discussion following Proposition 3.3, in \cite{Q}.
\item Here we are looking at congruences between modular forms of the same weight (i.e. $2$), and how to transfer them to half-integral weight.
For work on the analogous question for congruences between forms of different weights, see \cite{D} (which uses work of Stevens \cite{St} to go beyond special cases), and \cite[Theorem 1.4]{MO} for a different approach by McGraw and Ono.
\item The formula for $\WWW(\phi)$ used by B\"ocherer and Schulze-Pillot has coefficient of $\theta_i$ equal to $\frac{\phi(y_i)}{w_i}$ rather than just $\phi(y_i)$, and their $\langle\phi,\phi\rangle$ has $w_i$ in the denominator (as in \cite[(6.2)]{G2}) rather than in the numerator. This is because our $\phi(y_i)$ is the same as their $\phi(y_i)/w_i$. Their $\phi$ is an eigenvector for standard Hecke operators $T_p$ defined using right translation by double cosets (as in \cite[(6.6)]{G2}), which are represented by Brandt matrices, and are self-adjoint for their inner product. The Hecke operators we use below are represented by the transposes of Brandt matrices (as in \cite[Proposition 4.4]{G2}), and are self-adjoint for the inner product we use here (see the final remark). This accounts for the adjustment in the eigenvectors.
\end{enumerate}

\section{Proof of Theorem \ref{main} }
First, from $a_p(f)\equiv a_p(g)\pmod{\lambda}~~~\forall\text{ primes } p\nmid N\ell,$ it follows that $a_p(f)\equiv a_p(g)\pmod{\lambda}~~~\forall\text{ primes } p.$ This is because the congruence for $p\nmid N\ell$ (or even just for almost all primes) gives an ismorphism of residual global Galois representations. The Hecke eigenvalues at the remaining primes $p$ can be recovered from the Galois representations, from the unramified quotient of the restriction to $\Gal(\Qbar_p/\QQ_p)$ when $p\neq \ell$. (For $p=\ell$ this also applies in the ordinary case, by a theorem of Deligne \cite[Theorem 2.5]{Ed}, and in the supersingular case $a_{\ell}(f)\equiv a_{\ell}(g)\equiv 0\pmod{\lambda}$.)

Let $\TT$ be the $\ZZ$-algebra generated by the linear operators $T_p$ (for primes $p\nmid N$) and $U_p$ (for primes $p\mid N$) on the $q$-new subspace $S_2(\Gamma_0(N))^
{q\text{-}
\mathrm{new}}$ (the orthogonal complement of the subspace of those old forms coming from $S_2(\Gamma_0(N/q))$). Let $K$ be a number field sufficiently large to accommodate all the Hecke eigenvalues $a_p(f)$ and $a_p(g)$. The homomorphism $\theta_f:\TT\rightarrow K$ such that $T_p\mapsto a_p(f)$ and $U_p\mapsto a_p(f)$ has kernel $\pp_f$, say, and similarly we can define $\theta_g$ and $\pp_g$. Let $\lambda$ be a prime ideal of $O_K$, dividing a rational prime $\ell>3$ with $\ell\nmid N(q-1)$, such that $a_p(f)\equiv a_p(g)\pmod{\lambda}$ for all $p$. The homomorphism $\overline{\theta_f}=\overline{\theta_g}:\TT\rightarrow \FF_{\lambda}:=O_K/\lambda$ such that $\overline{\theta_f}(t)=\overline{\theta_f(t)}$ for all $t\in \TT$, has a kernel $\mm$ which is a maximal ideal of $\TT$, containing $\pp_f$ and $\pp_g$, with $k_{\mm}:=\TT/\mm\subseteq \FF_{\lambda}$.

By \cite[Theorem 3.10]{R}, $\TT$ may be viewed as a ring of endomorphisms of the toric part $T$ of the special fibre at $q$ of the N\'eron model of the Jacobian of the modular curve $X_0(N)$. Hence $\TT$ may be viewed as a ring of endomorphisms of the character group $X:=\Hom(T,\mathbb{G}_m)$, which is naturally identified with the set of divisors of degree zero (i.e $\ZZ$-valued functions summing to $0$) on the finite set $B^{\times}(\QQ)\backslash B^{\times}(\A_f)/\hat{R}$. (This finite set is in natural bijection with the set of intersection points of the two irreducible components of the special fibre at $q$ of a regular model of $X_0(N)$, with $R_i$ isomorphic to the endomorphism ring of the corresponding supersingular elliptic curve.) This gives a geometrical realisation of the Jacquet-Langlands correspondence. We find eigenvectors $\phi_f$ and $\phi_g$, on which $\TT$ acts through $\TT/\pp_f$ and $\TT/\pp_g$ respectively, in $X\otimes_{\ZZ}K$. We may extend coefficients to $K_{\lambda}$, and scale $\phi_f$ and $\phi_g$ to lie in $X\otimes O_{\lambda}$ but not in $\lambda (X\otimes O_{\lambda})$.

By \cite[Theorem 6.4]{R} (which uses the conditions that $\rhobar_{f,\lambda}$ is irreducible and that $\ell\nmid N(q-1)$), $\dim_{k_{\mm}}(X/\mm X)\leq 1$. The proof of this theorem of Ribet uses his generalisation to non-prime level \cite[Theorem 5.2(b)]{R} of Mazur's ``multiplicity one'' theorem that $\dim_{k_{\mm}}(J_0(N)[\mm])=2$ \cite[Proposition 14.2]{M}, and Mazur's level-lowering argument for $q\not\equiv 1\pmod{\ell}$. We can remove the condition $\ell\nmid N$ in the case $\ell=q$, using Wiles's further generalisation of Mazur's multiplicity one theorem \cite[Theorem 2.1(ii)]{Wi}.

We can localise at $\mm$ first, so $\phi_f,\phi_g\in X_{\mm}\otimes O_{\lambda}$ and $\dim_{k_{\mm}}(X_{\mm}/\mm X_{\mm})\leq 1$.
In fact, since we are looking only at a Hecke ring acting on $q$-new forms (what Ribet calls $\TT_1$), we must have $\dim_{k_{\mm}}(X_{\mm}/\mm X_{\mm})=1$. It follows from \cite[Theorem 3.10]{R}, and its proof, that $X_{\mm}\otimes_{\ZZ_{\ell}}\QQ_{\ell}$ is a free $\TT_{\mm}\otimes_{\ZZ_{\ell}}\QQ_{\ell}$-module of rank $1$. Then an application of Nakayama's Lemma shows that $X_{\mm}$ is a free $\TT_{\mm}$-module of rank $1$. Now $\TT_{\mm}$ is a Gorenstein ring, as in \cite[Corollary 15.2]{M}, so $\dim_{k_{\mm}}((\TT_{\mm}/\ell \TT_{\mm})[\mm])=1$ (by \cite[Proposition 1.4(iii)]{T}) and hence $\dim_{k_{\mm}}((X_{\mm}/\ell X_{\mm})[\mm])=1$. It follows by basic linear algebra that $((X_{\mm}\otimes_{\ZZ_{\ell}}O_{\lambda})/\ell (X_{\mm}\otimes_{\ZZ_{\ell}}O_{\lambda}))[\mm\otimes_{\ZZ_{\ell}}O_{\lambda}]$ is a free $(k_{\mm}\otimes_{\FF_{\ell}}\FF_{\lambda})$-module of rank $1$, using the assumption that $K_{\lambda}/\QQ_{\ell}$ is unramified.

Now $(k_{\mm}\otimes_{\FF_{\ell}}\FF_{\lambda})\simeq \prod_{k_{\mm}\hookrightarrow\FF_{\lambda}}\FF_{\lambda}$, and it acts on both $\phi_f$ and $\phi_g$ through the single component corresponding to the map $k_{\mm}\hookrightarrow\FF_{\lambda}$ induced by $\overline{\theta_f}=\overline{\theta_g}$. Hence $\phi_f$ and $\phi_g$ reduce to the same $1$-dimensional $\FF_{\lambda}$-subspace of $(X_{\mm}\otimes_{\ZZ_{\ell}}O_{\lambda})/\ell (X_{\mm}\otimes_{\ZZ_{\ell}}O_{\lambda})$, and by rescaling by a $\lambda$-adic unit, we may suppose that their reductions are the same, i.e. that $\phi_f(y_i)\equiv\phi_g(y_i)\pmod{\lambda}\,\,\,\forall\,\,i$. The claimed congruence is now an immediate consequence of $\WWW(\phi)=\sum_{i=1}^h \phi(y_i)\theta_i,$ and the integrality of the Fourier coefficients of the $\theta_i$.

\section{Two examples}
$\mathbf{N=170}.$ Let $f$ and $g$ be the newforms for $\Gamma_0(170)$ attached to the isogeny classes of elliptic curves over $\QQ$ labelled $\mathbf{170b}$ and $\mathbf{170e}$ respectively, in Cremona's data \cite{C}. For both $f$ and $g$ the Atkin-Lehner eigenvalues are $w_2=w_5=+1$, $w_{17}=-1$. The modular degrees of the optimal curves in the isogeny classes $\mathbf{170b}$ and $\mathbf{170e}$ are $160$ and $20$, respectively. Both are divisible by $5$, with the consequence that $5$ is a congruence prime for $f$ in $S_2(\Gamma_0(170))$, and likewise for $g$. In fact $f$ and $g$ are congruent to each other mod $5$.
\vskip10pt
\begin{tabular}{|c|c|c|c|c|c|c|c|c|c|c|c|c|c|}\hline $p$ & $3$ & $7$ & $11$ & $13$ & $19$ & $23$ & $29$ & $31$ & $37$ & $41$ & $43$ & $47$ & $53$\\\hline $a_p(f)$ & $-2$ & $2$ & $6$ & $2$ & $8$ & $-6$ & $-6$ & $2$ & $2$ & $-6$ & $-4$ & $12$ & $6$\\\hline
$a_p(g)$ & $3$ & $2$ & $-4$ & $-3$ & $3$ & $-6$ & $9$ & $-3$ & $-8$ & $-6$ & $6$ & $-13$ & $-9$\\\hline\end{tabular}

The Sturm bound \cite{Stu} is $\frac{kN}{12}\prod_{p\mid N}\left(1+\frac{1}{p}\right)=54$, so the entries in the table (together with the Atkin-Lehner eigenvalues) are sufficient to prove the congruence $a_n(f)\equiv a_n(g)$ for all $n\geq 1$.

Using the computer package Magma, one can find matrices for Hecke operators acting on the Brandt module for $D=17$, $M=10$, for which $h=24$. Knowing in advance the Hecke eigenvalues, and computing the null spaces of appropriate matrices, one easily finds that we can take $[\phi_f(y_1),\ldots,\phi_f(y_{24})]$ and $[\phi_g(y_1),\ldots,\phi_g(y_{24})]$ (with the ordering as given by Magma) to be $$[-4,-4,-4,-4,5,5,5,5,5,5,5,5,2,2,-1,-1,-1,-1,-1,-1,-1,-1,-10,-10]$$ and $$[1,1,1,1,0,0,0,0,0,0,0,0,2,2,-1,-1,-1,-1,-1,-1,-1,-1,0,0]$$ respectively, and we can observe directly a mod $5$ congruence between $\phi_f$ and $\phi_g$.

Using the computer package Sage, and Hamieh's function \newline
``shimura\_lift\_in\_kohnen\_subspace'' \cite[\S 4]{H}, we found
$$\WWW(\phi_f)=-4q^{20}+16q^{24}-24q^{31}+16q^{39}+20q^{40}+8q^{56}-8q^{71}-40q^{79}+4q^{80}+16q^{95}$$ $$-16q^{96}+O(q^{100}),$$
$$\WWW(\phi_g)=-4q^{20}-4q^{24}-4q^{31}-4q^{39}+8q^{56}+12q^{71}+4q^{80}-4q^{95}+4q^{96}+O(q^{100}),$$
in which the mod $5$ congruence is evident. Unfortunately the condition $\ell\nmid N$ (or $\ell=D$) does not apply in this example.

$\mathbf{N=174}.$ Let $f$ and $g$ be the newforms for $\Gamma_0(174)$ attached to the isogeny classes of elliptic curves over $\QQ$ labelled $\mathbf{174a}$ and $\mathbf{174d}$ respectively, in Cremona's data \cite{C}. For both $f$ and $g$ the Atkin-Lehner eigenvalues are $w_2=w_{29}=+1$, $w_{3}=-1$. The modular degrees of the optimal curves in the isogeny classes $\mathbf{174a}$ and $\mathbf{174d}$ are $1540$ and $10$, respectively. Both are divisible by $5$, with the consequence that $5$ is a congruence prime for $f$ in $S_2(\Gamma_0(174))$, and likewise for $g$. In fact $f$ and $g$ are congruent to each other mod $5$.
\vskip10pt
\begin{tabular}{|c|c|c|c|c|c|c|c|c|c|c|c|c|c|c|}\hline $p$ & $5$ & $7$ & $11$ & $13$ & $17$ & $19$ & $23$ & $31$ & $37$ & $41$ & $43$ & $47$ & $53$ & $59$\\\hline $a_p(f)$ & $-3$ & $5$ & $6$ & $-4$ & $3$ & $-1$ & $0$ & $-4$ & $-1$ & $-9$ & $-7$ & $-3$ & $-6$ & $3$\\\hline
$a_p(g)$ & $2$ & $0$ & $-4$ & $6$ & $-2$ & $4$ & $0$ & $-4$ & $-6$ & $6$ & $-12$ & $-8$ & $-6$ & $8$\\\hline\end{tabular}

The Sturm bound \cite{Stu} is $\frac{kN}{12}\prod_{p\mid N}\left(1+\frac{1}{p}\right)=60$, so the entries in the table (together with the Atkin-Lehner eigenvalues) are sufficient to prove the congruence $a_n(f)\equiv a_n(g)$ for all $n\geq 1$.

Using Magma, one can find matrices for Hecke operators acting on the Brandt module for $D=3$, $M=58$, for which $h=16$. We find $[\phi_f(y_1),\ldots,\phi_f(y_{16})]$ and $[\phi_g(y_1),\ldots,\phi_g(y_{16})]$ to be $$[2,2,-5,-5,-5,-5,10,10,10,10,-2,-2,-2,-2,-8,-8]$$ and $$[2,2,0,0,0,0,0,0,0,0,-2,-2,-2,-2,2,2]$$ respectively, and we can observe directly the mod $5$ congruence between $\phi_f$ and $\phi_g$ proved on the way to Theorem \ref{main}.

Using the computer package Sage, and Hamieh's function \newline
``shimura\_lift\_in\_kohnen\_subspace'' \cite[\S 4]{H}, we found (with appropriate scaling)
$$\WWW(\phi_f)=2q^4-10q^7-2q^{16}-8q^{24}+10q^{28}+2q^{36}+20q^{52}-10q^{63}+2q^{64}-12q^{87}-4q^{88}$$ $$+8q^{96}+O(q^{100}),$$
$$\WWW(\phi_g)=2q^4-2q^{16}+2q^{24}+2q^{36}+2q^{64}-2q^{87}-4q^{88}-2q^{96}+O(q^{100}),$$
in which the mod $5$ congruence is evident.
The condition $\ell\nmid N(q-1)$ does apply to this example, and $\rhobar_{f,\ell}$ is irreducible, since we do not have $a_p(f)\equiv 1+p\pmod{\ell}$ for all $p\nmid \ell N$.

{\em Remark. } The norm we used comes from a bilinear pairing $\langle,\rangle:X\times X\rightarrow \ZZ$ such that $\langle y_i,y_j\rangle=w_j\delta_{ij}$. The Hecke operators $T_p$ for $p\nmid N$ are self-adjoint for $\langle,\rangle$, since if $E_i$ is the supersingular elliptic curve associated to the class representated by $y_i$, then $\langle T_p y_i,y_j\rangle$ is the number of cyclic $p$-isogenies from $E_i$ to $E_j$, while $\langle y_i, T_p y_j\rangle$ is the number of cyclic $p$-isogenies from $E_j$ to $E_i$, and the dual isogeny shows that these two numbers are the same. See the discussion preceding \cite[Proposition 3.7]{R}, and note that the factor $w_j=\#\Aut(E_j)$ intervenes between counting isogenies and just counting their kernels.

We have $\phi_f-\phi_g=\lambda\phi$ for some $\phi\in X\otimes O_{\lambda}$. Hence $\phi=\frac{1}{\lambda}(\phi_f-\phi_g)$. Now $\phi_f$ and $\phi_g$ are simultaneous eigenvectors for all the $T_p$ with $p\nmid N$, and are orthogonal to each other, so we must have $\frac{1}{\lambda}=\frac{\langle\phi,\phi_f\rangle}{\langle\phi_f,\phi_f\rangle}$. Consequently, $\lambda\mid\langle\phi_f,\phi_f\rangle$, and similarly $\lambda\mid\langle\phi_g,\phi_g\rangle$. We can see this directly in the above examples, where $\lambda=\ell=5$. In the first one, the GramMatrix command in Magma shows that all $w_i=2$, so $\langle\phi_f,\phi_f\rangle=960$ and $\langle\phi_g,\phi_g\rangle=40$. In the second example, $w_1=w_2=4$ while
$w_i=2$ for all $3\leq i\leq 16$, so $\langle\phi_f,\phi_f\rangle=1320$ and $\langle\phi_g,\phi_g\rangle=80$.

{\em Acknowledgments. } We are grateful to A. Hamieh for supplying the code for her Sage function
``shimura\_lift\_in\_kohnen\_subspace'', as used in \cite[\S 4]{H}. The second named author would like to thank
Sheffield University for the great hospitality.
She was supported by a DST-INSPIRE grant.

\end{document}